\documentclass{amsart}

\usepackage{amssymb}
\usepackage{eepic}
\usepackage{color}
\usepackage{mathtools,xparse}
\usepackage[linktocpage=true]{hyperref}
\usepackage[utf8x]{inputenc}

\allowdisplaybreaks

\numberwithin{equation}{section}

\newtheorem{theorem}{Theorem}[section]
\newtheorem{lemma}[theorem]{Lemma}

\newtheorem{corollary}[theorem]{Corollary}
\newtheorem{remark}[theorem]{Remark}

\newtheorem{TheoA}{Theorem A}

\newtheorem{TheoBB}{Theorem B1}
\newtheorem{TheoBBB}{Theorem B2}

\newcommand{\Z}{\mathbf{Z}}
\newcommand{\R}{\mathbf{R}}
\newcommand{\T}{\mathbf{T}}
\newcommand{\C}{\mathbf{C}}

\newcommand{\B}{\mathcal{B}}
\newcommand{\SL}{S \hskip-1pt L_n(\R)}

\def\G{\mathrm{G}}
\def\la{\langle}
\def\ra{\rangle}
\def\ii{\mathrm{i}}
\def\1{\mathbf{1}}
\def\H{\mathcal{H}}

\def\M{\mathcal{M}}

\newcommand{\dem}{\noindent {\bf Proof. }}

\newcommand{\demA}{\noindent {\bf Proof of Theorem A. }}
\newcommand{\demB}{\noindent {\bf Proof of Theorem B1 i). }}
\newcommand{\demBB}{\noindent {\bf Proof of Theorem B1 ii). }}
\newcommand{\demBBB}{\noindent {\bf Proof of Theorem B2 ii). }}
\newcommand{\demBBBB}{\noindent {\bf Proof of Theorem B2 i). }}

\newcommand{\fin}{\hspace*{\fill} $\square$ \vskip0.2cm}

\begin{document}

\null

\vskip20pt

\null

\begin{center}
{\huge Riesz-Schur transforms}

\vskip15pt

{\sc {Adri\'an M. Gonz\'alez-P\'erez, Javier Parcet \\ Jorge P\'erez Garc\'ia and \'Eric Ricard}}
\end{center}

\title[Riesz-Schur transforms]{}


\maketitle

\null

\vskip-45pt

\null

\begin{center}
{\large {\bf Abstract}}
\end{center}

\vskip-25pt

\null

\begin{abstract}
We investigate nontrigonometric forms of Riesz transforms in the context of Schur multipliers. This refines Grothendieck-Haagerup's endpoint criterion with a new condition for the Schatten $p$-boundedness of Schur multipliers and strengthens Potapov/Sukochev's solution of Arazy's conjecture. We recover as well dimension-free estimates for trigonometric Riesz transforms. Our discrete approach is much simpler than previous harmonic analysis and probabilistic approaches. As an application, we find a very simple proof of recent criteria for Schur multipliers of H\"ormander-Mikhlin and Marcinkiewicz type.  
\end{abstract}

\addtolength{\parskip}{+1ex}

\vskip25pt

\section*{\bf \large Introduction}

Schur multipliers are fundamental maps on matrix algebras, which admit a rather simple definition $S_M(A) = (M(j,k) A_{jk})_{jk}$ for any matrix $A \in \B(\ell_2(\Z))$ and certain symbol $M: \Z \times \Z \to \C$. This easily extends to nonatomic $\sigma$-finite measure spaces $(\Omega,\mu)$, by restricting to operators $A$ in $L_2(\Omega,\mu)$ admitting a kernel representation over $\Omega \times \Omega$. Schur multipliers have shown deep connections with functional and harmonic analysis, operator algebras, and geometric group theory. They played a key role in Grothendieck's work on Banach spaces \cite{G,PisSim,PisBAMS}, Haagerup's invariant approximation methods for the group von Neumann algebras of free groups and rank one semisimple lattices \cite{DCH,H,CH}, Lafforgue/de la Salle's rigidity theorem for high-rank lattices \cite{dLdlS,LdlS} or Potapov/Sukochev's solution to Krein's conjecture on operator-Lipschitz functions by solving the stronger Arazy's conjecture on Schur multipliers \cite{PS}. More recently, Schur multipliers have been instrumental in \cite{PRS} to construct singular Fourier multipliers in the group algebra of $\SL$ and this has motivated a renewed interest in Schatten $p$-estimates for nonToeplitz Schur multipliers. Namely, this class is more stable under cut/paste or deformation of symbols than Fourier multipliers and yields important nontrigonometric extensions of fundamental results for them. This includes the nonToeplitz forms of celebrated multiplier theorems by Marcinkiewicz, H\"ormander-Mikhlin and Fefferman \cite{CLM,CGPT1,PST} and strong applications have been found in Lie group algebras \cite{CGPT2,PST}. Our results below extend these (best-known to date) $S_p$-estimates for Schur multipliers. 

\vskip10pt

\hrule {\footnotesize Adri\'an Gonz\'alez-P\'erez was partially supported by Ramón y Cajal grant RYC2022-037045-I (Ministerio de Ciencia, Spain). Adri\'an Gonz\'alez-P\'erez, Javier Parcet and Jorge P\'erez-Garc\'ia were supported by the Spanish Grant PID2022-141354NB-I00 and the Severo Ochoa Grant CEX2023-001347-S. \'Eric Ricard was supported by the French Grant ANR-19-CE40-0002.} 
  
It is important to recall a key relation between Fourier and Schur multipliers over amenable groups. Assume that $\G$ is an amenable group equipped with its Haar measure $\mu$. Let $S_p(L_2(\G,\mu))$ denote the Schatten $p$-class over $L_2(\G,\mu)$ and let $L_p(\mathcal{L}(\G))$ be the noncommutative $L_p$-space associated to the group von Neumann algebra $\mathcal{L}(\G)$. Let $M(g,h) = m(gh^{-1})$ be a Herz-Schur (also known as Toeplitz) symbol. If $e_{g,h}$ and $\lambda(g)$ respectively denote matrix units and the left regular representation over $\G$, consider the Schur multiplier $S_M: e_{g,h} \mapsto M(g,h) e_{g,h}$ and the Fourier multiplier $T_m: \lambda(g) \mapsto m(g) \lambda(g)$. Then, we know from \cite{CdlS,NR} that the following identity holds for $1 \le p \le \infty$ 
\begin{equation} \tag{FS} \label{Eq-FS}
\big\| S_M \hskip-2pt : S_p(L_2(\G,\mu)) \to S_p(L_2(\G,\mu)) \big\|_{\mathrm{cb}} = \big\| T_m \hskip-2pt : L_p(\mathcal{L}(\G)) \to L_p(\mathcal{L}(\G)) \big\|_{\mathrm{cb}}. 
\end{equation}
This is known as Fourier-Schur transference. It clarifies why Schatten $p$-bounds for nonToeplitz Schur multipliers over the classical groups $\R^n$ and $\T^n$ appear as nontrigonometric extensions of Fourier $L_p$-multipliers.

Riesz transforms are the archetypes of singular integrals, showing a great impact in harmonic analysis, fluid mechanics, differential geometry, and geometric measure theory. Given a square-integrable function $f: \R^n \to \C$, the Euclidean Riesz transforms acting on $f$ are defined as follows $$R_jf(x) = \partial_j (-\Delta)^{-\frac12}f(x) = \int_{\R^n} \Big\langle \frac{\xi}{|\xi|}, e_j \Big\rangle \widehat{f}(\xi) e^{2\pi i \langle x,\xi \rangle} \, d\xi.$$ Dimension-free estimates for the full Riesz transform $\mathcal{R}f = \nabla (-\Delta)^{-\frac12}f$ go back to the work of Gundy-Varopoulos and Stein \cite{GV,S}. Shortly after, a probabilistic interpretation by P.A. Meyer pioneered significant generalizations \cite{B1,B2,Gu,M,P0} for diffusion processes. The remarkable work of Lust-Piquard incorporated abelian groups and the first noncommutative models to the picture|her analysis for the Walsh system was applied by Naor in a striking connection with metric geometry \cite{LP1,LP2,LP3,N}. Riesz transforms in von Neumann algebras have been systematically investigated in \cite{JM,JMP}, the latter includes dimension-free estimates for arbitrary Markov processes in group von Neumann algebras. Given an index set $\Gamma$ and a Hilbert space $\H$, in this paper we consider Riesz-Schur transforms defined for finitely supported matrices $A \in \B(\ell_2(\Gamma))$ as follows $$\mathcal{R}_u(A) = \Big( \frac{u_j-u_k}{\|u_j-u_k\|} A_{jk} \Big)_{j,k \in \Gamma} \quad \mbox{for some} \quad u = (u_j)_{j \in \Gamma} \subset \H.$$ In what follows $S_p(\Gamma)$ denotes the Schatten $p$-class over $\ell_2(\Gamma)$. We start with dimension-free estimates for Riesz-Schur transforms. Instead of the usual statement in discrete settings \cite{JMP}, we allow non-mean-zero elements adding a diagonal term. 

\begin{TheoA}
Let $\Gamma$ be any index set and fix any family $(u_j)_{j \in \Gamma}$ of vectors in some Hilbert space $\H$. Consider a semifinite von Neumann algebra $(\M,\tau)$. Then the following inequalities hold for $x = \sum_{j,k} x_{jk} \otimes e_{jk} \in L_p(\M; S_p(\Gamma))$ and any $1 < p < \infty$ up to universal constants  
\begin{equation} \label{EqRiesz1} \tag{RS$_1$}
\Big\| \hskip-2pt \sum_{j,k \in \Gamma} x_{jk} \otimes e_{jk} \otimes \frac{u_j - u_k}{\|u_j - u_k\|} \Big\|_{RC_p} \hskip-1.5pt \lesssim c_p \Big\| \hskip-2pt \sum_{j,k \in \Gamma} x_{j,k} \otimes e_{j,k} \Big\|_p,
\end{equation}
\vskip-10pt 
\begin{equation} \label{EqRiesz2} \tag{RS$_2$} \Big\| \hskip-2pt \sum_{j,k \in \Gamma} x_{jk} \otimes e_{jk} \Big\|_p \lesssim \Big\| \hskip-4pt \sum_{\begin{subarray}{c} j,k \in \Gamma \\ u_j = u_k \end{subarray}} x_{jk} \otimes e_{jk} \Big\|_p + c_{p'} \Big\| \hskip-2pt \sum_{j,k \in \Gamma} x_{jk} \otimes e_{jk} \otimes \frac{u_j - u_k}{\|u_j - u_k\|} \Big\|_{RC_p}.
\end{equation} 
\vskip-8pt
\noindent Here $c_p = \max \big\{p, (\frac{p}{p-1})^{\frac32} \big\}$ and we adopt the convention that $\frac00 = 0$ when $u_j = u_k$.
\end{TheoA}

The $RC_p$-norm is nothing but a matrix analogue of the $L_p$-norm of a square function. By well-known approximation properties of Schur multipliers \cite{LdlS}, our discrete statement over $\ell_2(\Gamma)$ holds as well for general $\sigma$-finite spaces. Theorem A is a nonToeplitz/basis-free extension of the main result in \cite{JMP}. Indeed, it recovers the dimension-free estimates from \cite{JMP,LP3,S} by Fourier-Schur transference results \cite{CdlS,NR,PRS}. The Euclidean case follows by taking $\H = \R^n$ and $u_x=x$, followed by Theorem 1.19 of \cite{LdlS} and Fourier-Schur transference \eqref{Eq-FS}. Other $L_p$-estimates for different laplacians and/or groups follow using appropriate cocycle maps. As it is customary, the extra \lq diagonal\rq${}$ term in \eqref{EqRiesz2} disappears in continuous settings. Our constants match the best-known so far in this context. Arhancet and Kriegler also investigated it in \cite[Theorem 3.3]{AK}. They followed \cite{JMP} estimating $L_p$-norms of the gaussian gradient $\delta_\psi$ |equivalently the Dirac operator| and their statement is equivalent to ours via a Khintchine type inequality. The main contribution in Theorem A is a cleaner statement and prominently a much simpler proof, which has facilitated to uncover interesting applications below. None of the usual harmonic analysis or probabilistic methods |like Fourier transforms and Calder\'on-Zygmund techniques or diffusion/Markov semigroups and Pisier's reduction formula| are needed. On the contrary, our argument is modeled on Grothendieck's inequality \cite[Chapter 5]{PisSim}. This becomes quite enlightening and it is not accidental that Theorem A leads to a condition for Schatten $p$-boundedness of Schur multipliers refining the Grothendieck-Haagerup's endpoint criterion. 

\begin{TheoBB}
Consider arbitrary vectors $\{u_j, u_j', w_j, w_j': j \in \Gamma\}$ in some Hilbert space $\H$ and let $\Lambda: \H \to \H$ be a contraction. Then, the following inequalities hold for any $1 < p < \infty$ up to absolute constants$\hskip1pt :$  
\begin{itemize}
\item[i)] \emph{Grothendieck-Haagerup type $S_p$-criterion.} The symbols $$M(j,k) = \Big\langle \frac{u_j + u_k'}{\|u_j + u_k'\|}, \Lambda \Big( \frac{w_j + w_k'}{\|w_j + w_k'\|}\Big)\Big\rangle$$ define cb-bounded Schur multipliers with $\| S_M \|_{\mathrm{cb}(S_p(\Gamma))} \lesssim \max \big\{p,\frac{p}{p-1} \big\}^{\frac52}$.

\vskip5pt

\item[ii)] \emph{Square root of Arazy's divided differences.} Given a nondecreasing Lipschitz function $f: \R \to \R$, it turns out that the symbol $$M_f(x,y) = \sqrt{\frac{f(x) - f(y)}{x-y}}$$ defines a cb-bounded multiplier with $\| S_{M_f} \|_{\mathrm{cb}(S_p(\R))} \lesssim \max \big\{ p, \frac{p}{p-1} \big\}^{\frac52} \|f\|_{\mathrm{lip}}^{\frac12}$.
\end{itemize}
\end{TheoBB}

Taking $u_k' = w_j = 0$, we obtain symbols of the form $\langle \xi_j, \psi_k \rangle$ for some uniformly bounded families of vectors in $\H$. These characterize $S_\infty$-bounded Schur multipliers \cite{G,PisSim,PisBAMS} and Theorem B1 i) gives a weaker form of the Grothendieck-Haagerup criterion for Schatten $p$-classes. We also get better constants for certain particular cases. The second criterion yields a significant strengthening of Arazy's conjecture \cite{A,PS}. It should be compared with \cite[Corollary 3.4]{CGPT1} for $\alpha$-divided differences over $\alpha$-H\"older functions. Both results suggest in different forms that $S_p$-boundedness holds for $1 < p< \infty$ under weak forms of Lipschitz continuity. More details can be found in Remark \ref{Remalphabeta}. Our last result provides a very simple proof of two recent criteria which have received a lot of attention.    

\begin{TheoBBB}
Given $1 < p < \infty \hskip-2pt :$  
\begin{itemize}
\item[i)] \emph{Marcinkiewicz-Schur multipliers}. If $M:\Z \times \Z \to \C$ $$\big\| S_M \big\|_{\mathrm{cb}(S_p(\mathbf{Z}))} \ \lesssim \ \sup_{j,k \in \mathbf{Z}} \mathbf{Var}_{|\cdot| \sim 2^k} \Big\{ \big|M(j + \cdot, j)\big| + \big| M(j, j + \cdot) \big| \Big\}.$$

\item[ii)] \emph{H\"ormander-Mikhlin-Schur multipliers.} If $M \in \mathcal{C}^{[\frac{n}{2}]+1}(\R^{2n} \setminus \{0\})$ $$\hskip24pt \big\| S_M \big\|_{\mathrm{cb}(S_p(\mathbf{R}^n))} \lesssim \sum_{|\gamma| \le [\frac{n}{2}] +1} \hskip-5pt \Big\| |x-y|^{|\gamma|} \Big\{ \big| \partial_x^\gamma M(x,y) \big| + \big| \partial_y^\gamma M(x,y) \big| \Big\} \Big\|_\infty.$$
\end{itemize}
\end{TheoBBB}

Our argument is again simpler and intrinsic, avoiding corepresentations and Euclidean harmonic analysis. In fact, we shall prove a stronger form of Theorem B2 ii) which strengthens the main result in \cite{CGPT1}. In conclusion, Riesz-Schur transforms recover and refine well-known estimates for Hilbert transforms |Theorem A for $\dim \H = 1$| and Riesz transforms, Grothendieck-Haagerup's multipliers, Arazy's divided differences, Marcinkiewicz-Schur and H\"ormander-Mikhlin-Schur multipliers \cite{CLM,CGPT1,G,JMP,LP3,PS,S}. To the best of our knowledge, these references give all known sufficient conditions on symbols for $S_p$-bounded Schur multipliers and $p \ge 1$.  

\section{\bf \large Hilbert-valued $L_p$-spaces} \label{S1}

Let $\Gamma$ be any index set and consider a semifinite von Neumann algebra $\M$ with a normal semifinite faithful trace $\tau$. Then, given $1 \le p \le \infty$, the corresponding noncommutative $L_p$-spaces are denoted by $L_p(\M,\tau)$ or $L_p(\M)$. In particular, the Schatten $p$-class $S_p(\Gamma)$ is nothing but $L_p(\B(\ell_2(\Gamma)), \mathrm{tr})$ and $L_p(\M \bar\otimes \B(\ell_2(\Gamma)), \tau \otimes \mathrm{tr})$ will be identified with a space of matrices $(x_{jk})_{j,k \in \Gamma}$ with coefficients in $L_p(\M)$ or with sums of the form $\sum_{j,k} x_{jk} \otimes e_{jk}$ where $e_{jk}$ are the standard matrix units. Any such matrix belongs to $L_p(\M \bar\otimes \B(\ell_2(\Gamma)), \tau \otimes \mathrm{tr})$ when $$\sup_{\begin{subarray}{c} \Sigma \subset \Gamma \\ \Sigma \ \mathrm{finite} \end{subarray}} \Big\| \sum_{j,k \in \Sigma} x_{jk} \otimes e_{jk} \Big\|_p < \infty.$$

\subsection{The spaces $L_p(\M, \H_{rc})$}

We recall basic facts on $p$-row/column operator spaces.  Our general reference is \cite{P}. We present them using a tensor
approach rather than a coordinate one. Let $\H$ be a complex Hilbert space whose scalar product is chosen to be antilinear on the left. Given $1\leq p\leq \infty$, we denote by $L_p(\M; \H_c)$ |or simply $C_p$ or $C_p(\M)$ if no confusion may occur| the 
associated $p$-column space with a fixed norm 1 vector $\xi \in \H$. Namely, if we identify every vector $u \in \H$ with the rank-one operator $\Lambda_u(h) = \langle \xi,h \rangle u$, then $L_p(\M; \H_c)$ consists of the closure in $L_p(\M \bar\otimes \B(\H))$ (weak-$*$ closure for $p=\infty$) of simple tensors. Given $\sum_{k=1}^n x_k\otimes u_k$ a simple tensor in $L_p(\M) \otimes \H$, its norm is
$$\Big\| \sum_{k=1}^n x_k \otimes u_k \Big\|_{C_p} = \, \Big\|\Big(\sum_{j,k=1}^n \la u_j, u_k\ra x_j^*x_k\Big)^{\frac12}\Big\|_p.$$ If $\Lambda :\H \to \H$ is a bounded map, then $\mathrm{id} \otimes \Lambda$ extends to a completely bounded map (weak-$*$ continuous for $p=\infty$) on $L_p(\M;\H_c)$ with the same norm. We say that $C_p$ is homogeneous. Similarly, the $p$-row Hilbert space $L_p(\M; \H_r) \subset L_p(\M\bar \otimes \B(\H))$ or simply $R_p$ consists of adjoints of elements in $C_p$. It is the closure (weak-$*$-closure for $p=\infty$) of simple tensors $\sum_{k=1}^n x_k\otimes \bar u_k \in L_p(\M) \otimes \overline{\H}$ identifying each $\bar u \in \overline{\H}$ with $\Lambda_u^*(h) = \la u,h\ra \xi$ in $\B(\H)$. The norm is then given by 
$$\Big\| \sum_{k=1}^n x_k \otimes \bar u_k \Big\|_{R_p} = \, \Big\| \Big( \sum_{j,k=1}^n \la u_j, u_k \ra x_jx_k^* \Big)^{\frac12} \Big\|_p.$$
This space is also homogeneous. If $\Lambda: \H \to \H$ is bounded and $\bar \Lambda(\bar h) := \overline{\Lambda(h)}$, then 
$\mathrm{id} \otimes \bar \Lambda$ extends to a completely bounded map (weak-$*$ continuous for $p=\infty$) on $L_p(\M; \H_r)$ with the same norm as $\Lambda$. Recall the duality $C_p^*=R_{p'}$ and $R_p^*=C_{p'}$ for $1\leq p<\infty$ and $p'=\frac p{p-1}$ under the bracket $$\big\la x \otimes u, y \otimes \bar w \big\ra_{C_p,R_{p'}} = \tau(xy) \la w, u\ra_\H.$$ 

To define the space $RC_p$, we need to see $(C_p,R_p)$ as a compatible couple of Banach spaces, which requires a linear identification between $\H$ and $\overline \H$. This boils down to the choice of an orthonormal basis $(e_j)_{j \in \Gamma}$ of $\H$ or, equivalently, a real subspace $\H_\R \subset \H$ such that $\H$ is a trivial complexification of $\H_\R$; meaning that $\H = \C \otimes_2 \H_\R$ as a real Hilbert space with the obvious $\C$-structure. Thus, there is an isometric conjugation $J: \H \to \H$ acting trivially on $\H_\R$. With this, we have continuous linear maps $C_p, R_p \to L_p(\M)^\Gamma$ respectively extending $$x \otimes u \mapsto \big(\la e_j,u\ra x\big)_{j \in \Gamma} \quad \mbox{and} \quad x \otimes \bar u \mapsto \big(\la e_j,J(u) \ra x\big)_{j \in \Gamma}.$$ Another way to say it is that $x \otimes u \in C_p$ corresponds to $x \otimes \overline{J(u)} \in R_p$. As usual we define $RC_p$ to be $R_p+C_p$ if $p\le2$ and $R_p \cap C_p$ if
$p\geq 2$. Thus, we find the following expressions for elementary tensors 
$$\Big\| \sum_{k=1}^n x_k \otimes u_k \Big\|_{RC_p} \!\! = \begin{cases} \displaystyle \inf \Big\{ \Big\| \sum_{k=1}^m y_k \otimes v_k \Big\|_{C_p} + \Big\| \sum_{k=1}^r z_k \otimes \overline{J(w_k)} \Big\|_{R_p} \Big\} \!\! & \!\! \mbox{for} \ 1 \le p \le 2, \\ \displaystyle \max \Big\{ \Big\| \sum_{k=1}^n x_k\otimes u_k \Big\|_{C_p}, \Big\| \sum_{k=1}^n x_k\otimes \overline{J(u_k)} \Big\|_{R_p} \Big\} \!\! & \!\! \mbox{for} \ 2 \le p \le \infty, \end{cases}$$
where the infimum runs over all decompositions $x= \sum_{k=1}^m y_k \otimes v_k + \sum_{k=1}^r z_k \otimes w_k$. It can easily be seen that one can indeed restrict to simple tensors. The space $RC_p$ is then the completion of simple tensors $L_p(\M) \otimes \H$ under these norms, we only consider $1<p<\infty$ here. $RC_p$ is again homogeneous, if $\Lambda: \H \to \H$ is a bounded map, then $\mathrm{id} \otimes \Lambda$ on simple tensors uniquely extends to a complete bounded map on $RC_p$ with cb-norm $\|\Lambda\|$. Because $J$ is isometric, this follows easily from the same facts for $R_p$, $C_p$ and the norm formulas. The duality between $RC_p^*=RC_{p'}$ for $1<p<\infty$ follows from that of row and column spaces up to identifications. The duality bracket is given for simple tensors by 
\begin{equation} \label{Eq-DualityRCp}
\big\la x \otimes u, y \otimes w \big\ra_{RC_p,RC_{p'}} = \tau(xy) \la {J(w)}, u \ra_\H.
\end{equation}

\subsection{Noncommutative Khintchine inequalities}

We shall be using gaussian formulations of Lust-Piquard's noncommutative Khintchine inequalities. Let us recall the gaussian or boson functor associated to a real Hilbert space  $\H_\R$ with trivial complexification $\H$. This is  a commutative von Neumann algebra $L_\infty(\Omega,\mu)$ equipped with a probability measure $\mu$ and generated by centered gaussian variables $\{W(h): h \in \H_\R\}$ with covariance $$\int_\Omega W(h_1) W(h_2) \, d\mu = \la h_1,h_2\ra$$ in the sense that $L_\infty(\Omega) = \{e^{\ii s W(h)}: h \in \H_\R, s \in \R \}''$. If $h \in \H$, then $h = a+\ii b$ with $a,b \in \H_\R$ and we set $W(h)=W(a)+\ii W(b)$. In what follows we shall use the operator-valued Khintchine inequalities in $L_p(\Omega)$. Here we consider the $RC_p$-spaces over $\H$ with $\H_\R$ as real particular subspace. The result below follows by a simple change of basis from the main result in \cite{LP0} and \cite[Remark III.7]{LPP}.

\begin{theorem}\label{khin}
There exists constants $A, B > 0$ such that for any semifinite von Neumann algebra $(\M,\tau)$, any real Hilbert space $\H_\R$, any $1< p<\infty$, any $n\geq 1$ and $x_k\in L_p(\M)$, $u_k \in \H$ for $k \leq n$ $$A^{-1}\Big\|\sum_{k=1}^n x_k\otimes u_k\Big\|_{RC_p} \leq \Big\|\sum_{k=1}^n x_k\otimes W(u_k)\Big\|_p \leq B \sqrt p \, \Big\|\sum_{k=1}^n x_k\otimes u_k\Big\|_{RC_p}.$$
\end{theorem}

Let $\mathrm{X}_p = \overline{\rm span} \{ x \otimes W(u): x \in L_p(\M), \, u \in \H \}\subset L_p(\M \bar\otimes L_\infty(\Omega))$.  Then, the Khintchine inequalities provide isomorphisms $\sigma_p : \mathrm{X}_p \to RC_p$ compatible with the duality bracket: for simple tensors $(f,g) \in L_p \times L_{p'}$ $$\big\la \sigma_p(f), \sigma_{p'}(g) \big\ra_{RC_p,RC_{p'}} = \int_\Omega \tau (f(w)g(w)) \, d\mu(w) \quad \mbox{for} \quad 1<p<\infty.$$ Moreover, taking the natural conditional expectation $$\mathbf{E}: L_p(\M \bar\otimes L_\infty(\Omega)) \ni f \mapsto \int_\Omega f(w) \, d\mu(w) \in L_p(\M),$$ we can relate $\mathrm{X}_p$ with $R_p$ and $C_p$ for $p > 2$ as follows $$\|x\|_{R_p} = \big\| \mathbf{E}(\sigma_p^{-1}(x) \sigma_p^{-1}(x)^*) \big\|_{p/2}^{1/2} \quad \mbox{and} \quad \|x\|_{C_p} = \big\| \mathbf{E} (\sigma_p^{-1}(x)^* \sigma_p^{-1}(x)) \big\|_{p/2}^{1/2}.$$ Let $Q$ be the gaussian projection. That is, the orthogonal projection of  $L_2(\Omega)$ onto the closed subspace $\{W(h): \, h \in \H\}$. For any von Neumann algebra $(\M,\tau)$, $\mathrm{id} \otimes Q$  is also well defined on $L_2(\M \bar\otimes L_\infty(\Omega))$ and it extends to an $L_p$-projection.
   
\begin{lemma}\label{Qproj}
The map $\mathrm{id} \otimes Q$ extends to a bounded projection on $L_p(\M \bar\otimes L_\infty(\Omega))$ for all $1<p<\infty$ with norm less than $c\sqrt{\max\{p,p/p-1\}}$ for some constant $c>0$.
\end{lemma}

\dem If $p>2$ and $x \in L_p\cap L_2$ is a simple tensor, then its gaussian projection $(\mathrm{id} \otimes Q)(x)$ is of the form $\sum_{k=1}^n x_k \otimes W(u_k)$ and we easily get the inequality below $$\big\|\sigma_p((\mathrm{id} \otimes Q)(x))\big\|_{C_p} = \big\|\mathbf{E} \big((\mathrm{id} \otimes Q)(x)^*(\mathrm{id} \otimes Q)(x)\big)\big\|_{p/2}^{1/2} \leq \big\|\mathbf{E} (x^*x)\big\|_{p/2}^{1/2}\leq \|x\|_p.$$ Similar inequalities hold for $R_p$ and we can conclude thanks to the Khintchine inequalities in Theorem \ref{khin} and density. The case $1<p<2$ follows by duality. \fin

\subsection{$ORC_p$-bounded families of operators}

Now we consider noncommutative forms of the notion of $R$-bounded family of operators \cite[Chapter 8]{HvNVW}. Let $(\delta_n)_{n\in \Z}$ stand for the canonical basis of $\ell_2(\Z)$ and fix $1\leq p\leq \infty$. Let $(\phi_n)_{n \in \Z}$ be a sequence of operators $\phi_n: L_p(\M) \to L_p(\M)$. We say that $(\phi_n)_{n \in \Z}$ is $OC_p$-bounded if there is a bounded operator $\phi$ on $C_p = L_p(\M;\ell_2(\Z)_c)$ such that $\phi (x \otimes \delta_n) = \phi_n(x) \otimes \delta_n$. A family $\Phi$ of operators $L_p(\M)\to L_p(\M)$ is $OC_p$-bounded if $(\phi_n)_{n\in \Z}$ is $OC_p$-bounded up to a uniformly bounded constant for all  sequences $(\phi_n)_{n\in \Z} \subset \Phi$. Similarly we may define $OR_p$ and $ORC_p$ families of operators. We shall give in Section \ref{S3.2} below a few examples related to Riesz-Schur transforms. The following result is classical, see Lemma 4.2 in \cite{JLMX}.

\begin{lemma}\label{absconv}
If a set $\Phi$ is $OC_p$-bounded, the closure of its absolute convex hull for the point-weak topology is $OC_p$-bounded with the same constant. The same holds for $OR_p$-bounded or $ORC_p$-bounded families of operators.
\end{lemma}

Let $(\Delta_n)_{n\in \Z}$ be a family of disjoint sets in $\Gamma \times \Gamma$. Assume that the sequence of Schur multipliers $(S_{\Delta_n})_{n \in \Z}$ with symbol $1_{(j,k) \in \Delta_n}$ is completely unconditional on $S_p$ for $1<p<\infty$ (in the strongest sense). Using the noncommutative Khintchine inequalities, we get for $x = (x_{jk})_{j,k} \in L_p(\M \bar \otimes \B(\ell_2(\Gamma)))$
\begin{equation} \label{LP}
\Big\| \sum_{n \in \Z} \sum_{j,k \in \Gamma} 1_{(j,k) \in \Delta_n} x_{jk} \otimes e_{jk} \otimes \delta_n \Big\|_{RC_p} \approx_{D_p} \Big\| \sum_{j,k \in \Gamma} x_{jk} \otimes e_{jk} \Big\|_p.
\end{equation}
The result below is straightforward and follows from the definitions above.

\begin{lemma}\label{RCbound}
Let $\Phi$ be an $RC_p$-bounded family of operators with constant $K_p$. Then for any sequence of Schur multipliers 
$(S_{M_n})_{n\in \Z}$ in $\Phi$ and any sequence $(\Delta_n)_{n\in \Z}$ satisfying \eqref{LP}, the symbol 
$$M(j,k) = \sum_{n \in \Z} 1_{(j,k) \in \Delta_n} M_n(j,k) \quad \mbox{for} \quad j,k \in \Gamma$$
defines a Schur multiplier on $S_p(\Gamma)$ with complete norm bounded by $K_p D_p$.
\end{lemma}

\section{\bf \large Square function inequalities}

In this section we prove the dimension-free estimates for Riesz-Schur transforms in Theorem A. As explained in the Introduction, our proof does not follow any of the standard approaches from harmonic analysis or probability theory. Then, we shall recall equivariant versions which show how to recover the trigonometric Riesz transform estimates from \cite{JMP} in the amenable case. 

\subsection{Riesz-Schur transforms}

We start with an easy transference argument. We use the usual ${\rm sgn}$ function (with ${\rm sgn}( 0)=1$) and fix the constant $C_p=\max \{p , p^{\prime}\}$.

\begin{lemma}\label{main0}
Let $(\Omega,\mu)$ be a probability measure space. Then, the following norm equivalence holds in $L_p(\M \bar\otimes \B(\ell_2(\Gamma)) \bar\otimes L_\infty(\Omega))$ for any $1<p<\infty$ and any \hskip-1pt family $(f_j)_{j \in \Gamma}$ of real-valued \hskip-1pt measurable \hskip-1pt functions on $\Omega$
$$
\Big\| \sum_{j,k \in \Gamma} x_{jk} \otimes e_{jk} \otimes \mathrm{sgn} {(f_j-f_k)} \Big\|_p \approx_{C_p} \Big\| \sum_{j,k \in \Gamma} x_{jk} \otimes e_{jk} \Big\|_p.
$$
\end{lemma}

\dem It suffices to assume $(x_{jk})$ is finitely supported. Since
$$\Big\| \sum_{j,k \in \Gamma} x_{jk} \otimes e_{jk} \otimes \mathrm{sgn} {(f_j-f_k)} \Big\|_p^p = \int_\Omega \Big\| \sum_{j,k \in \Gamma} \mathrm{sgn} {(f_j(\omega)-f_k(\omega))} x_{jk} \otimes e_{jk} \Big\|^p_p d\mu(\omega),$$ it suffices to check the pointwise equivalence for all $\omega\in \Omega$. This follows from the complete boundedness of the Hilbert transform as a Schur multiplier after a suitable row/column permutation. This completes the proof. \fin

\demA By approximation we can clearly assume that $\Gamma$ is finite, say $\Gamma = \{1,2, \ldots,n\}$ and let $M_n = \B(\ell_2(\Gamma))$. Assume for the moment that $u_j \in \H_\R$ and consider its gaussian functor. Since $W(u_j-u_k)$ is a gaussian variable with $L_2$-norm $\|u_j-u_k\|$ we have (assuming $u_j\neq u_k$) that $$Q\big( \mathrm{sgn} \big(W(u_j-u_k)\big)\big) = \frac 1 {\sqrt {2 \pi}} \int_\R {\rm sgn}(s) s e^{-s^2/2} ds \, \frac{W(u_j-u_k)}{\|u_j-u_k\|} = \sqrt{\frac 2 \pi} \frac{W(u_j-u_k)}{\|u_k-u_k\|}.$$ Let $x = \sum_{j,k} x_{jk} \otimes e_{jk}$ and consider 
\begin{eqnarray*}
\mathcal{R}(x) & = & \sum_{j,k \in \Gamma} x_{jk} \otimes e_{jk} \otimes \frac{u_j - u_k}{\|u_j - u_k\|}, \\ \H (x) & = & \sum_{j,k \in \Gamma} x_{jk} \otimes e_{jk} \otimes \mathrm{sgn} \big( W(u_j) - W(u_k) \big).
\end{eqnarray*} 
Then we use the identity $$\sigma_p^{-1}(\mathcal{R}(x)) = \sqrt{\frac \pi 2}(\mathrm{id} \otimes Q)(\H(x))$$ and Lemma \ref{main0} for $f_j = W(u_j)$ to deduce \eqref{EqRiesz1}. Indeed, we have 
\begin{eqnarray*}
\big\| \mathcal{R}(x) \big\|_{R_p} \!\!\!\! & = & \!\!\!\! \big\| \mathbf{E} \big( \sigma_p^{-1}(\mathcal{R}(x)) \sigma_p^{-1}(\mathcal{R}(x))^* \big) \big\|_{p/2}^{1/2} \\ \!\!\!\! & = & \!\!\!\! \sqrt{\frac \pi 2} \big\| \mathbf{E} \big( (\mathrm{id} \otimes Q)(\H(x)) (\mathrm{id} \otimes Q)(\H(x))^* \big) \big\|_{p/2}^{1/2} \lesssim \| \H(x) \|_{p} \lesssim C_p \|x\|_p.
\end{eqnarray*}   
The same holds for the $C_p$-norm and \eqref{EqRiesz1} follows for $p \ge 2$. When $1< p<2$, we use the noncommutative Khintchine inequality in Theorem \ref{khin}, followed by Lemmas \ref{Qproj} and \ref{main0}. The use of Lemma \ref{Qproj} produces the asymmetry in the constants in the statement of Theorem A. The proof of \eqref{EqRiesz1} is complete for $u_j \in \H_\R$ and for the general case, let $\mathcal{K}_\R$ be $\H$ as a real Hilbert space with scalar product $\la \cdot,\cdot \ra_{\mathcal{K}_\R} = \mathrm{Re} \, \la \cdot,\cdot \ra_\H$. Let $\mathcal{K} = \mathcal{K}_\R + \ii \mathcal{K}_\R$ be the trivial complexification of $\mathcal{K}_\R$. Its norm is made so that for $u,w \in \mathcal{K}_\R$ we get $\|u+\ii w\|_\mathcal{K}^2=\|u\|_\H^2 +\|w\|_\H^2$. The map $\Lambda: \mathcal{K} \to \H$ given by $\Lambda(u+\ii w)=u+\ii w$ for $u,w \in \mathcal{K}_\R = \H$ is a well-defined $\C$-linear map with bound $\sqrt 2$. This gives \eqref{EqRiesz1} for $u_j \in \H$ (loosing a factor $\sqrt 2$ in the constant) since it holds for $u_j \in \mathcal{K}_\R$ and $RC_p$ is homogeneous.

Let us now consider the lower estimate \eqref{EqRiesz2}. If the matrix $(x_{jk})$ is supported by $\{(j,k) \in \Gamma \times \Gamma: u_j \neq u_k\}$, then the estimate follows by duality. Indeed, we know from \eqref{EqRiesz2} that the maps defined by $$\hskip1pt \mathcal{R}_p: L_p(\M \otimes M_n) \ni x \otimes e_{jk} \mapsto x \otimes e_{jk} \otimes \frac{u_j-u_k}{\|u_j-u_k\|} \in RC_p(\M \otimes M_n),$$
$$\overline{\mathcal{R}}_{p'}: L_{p'}(\M \otimes M_n) \ni x \otimes e_{jk} \mapsto x \otimes e_{jk} \otimes \frac{J(u_j-u_k)}{\|u_j-u_k\|} \in RC_p(\M \otimes M_n),$$ have norms $\lesssim c_p, c_{p'}$ respectively. Given a pair $(x,y) \in L_p(\M\otimes M_n) \times L_{p'}(\M\otimes M_n)$ and using the duality bracket \eqref{Eq-DualityRCp} for $RC_p$ and $RC_{p'}$, we deduce the following identity
$$\big\la \mathcal{R}_p(x),\overline{\mathcal{R}}_{p'}(y) \big\ra_{RC_p,RC_{p'}}=\sum_{j,k \in \Gamma} \tau(x_{jk}y_{kj}) \frac{\la u_j-u_k,u_k-u_j\ra}{\|u_j-u_k\|^2}=-({\rm tr}\otimes\tau)(xy),$$ with the support assumption. Taking the sup over all $y$ in the unit $p'$-ball gives $\|x\|_p \lesssim c_{p'} \|\mathcal{R}(x)\|_{RC_p}$. The general case follows by the triangular inequality. One can also note that the term $\sum_{u_j = u_k} x_{jk} \otimes e_{jk}$ is obtained from $\sum_{j,k} x_{jk} \otimes e_{jk}$ by taking a conditional expectation. This completes the proof. \fin

\subsection{Trigonometric Riesz transforms} We start with a simple remark which will allow us to connect our estimates for Riesz-Schur transforms in Theorem A to the dimension-free estimates for trigonometric Riesz transforms in \cite{JMP}. Recall that the constants $c_p$ in Theorem A are best-known so far.

\begin{lemma}\label{colmathom}
Let $p\geq 2$ and consider a family of contractions $\Lambda_j: \H \to \H$. Then the following inequality holds for $(x_{jk})_{j,k \in \Gamma} \in L_p(\M \otimes \B(\ell_2(\Gamma)))$ and any family of vectors $(\xi_{jk})_{j,k \in \Gamma}$ in $\H \hskip-2pt :$
$$\Big\| \sum_{j,k \in \Gamma} x_{jk} \otimes e_{jk} \otimes \Lambda_j(\xi_{jk}) \Big\|_{C_p} \le \Big\| \sum_{j,k \in \Gamma} x_{jk} \otimes e_{jk} \otimes \xi_{jk} \Big\|_{C_p}.$$ The same estimate holds in $R_p$ with contractions $\Lambda_k$ varying over columns $k \in \Gamma$.  
\end{lemma}

\dem This follows from the homogeneity and associativity of $C_p$. Just note that $e_{jk}$ can be expressed as a tensor $e_{j1} \otimes e_{1j}$, so that everything boils down to the boundedness of the diagonal map $e_{j1} \otimes u \mapsto e_{j1} \otimes \Lambda_j(u)$. The proof is complete. \fin

We are now in position to recover some results from \cite{JMP}. When $\Gamma = \G$ is a locally compact group, consider an orthogonal cocycle $\beta: \G \to \H_\R$ associated to an orthogonal action $\alpha$. That is, $\beta(gh) = \beta(g) + \alpha_g(\beta(h))$ and we find the  
\begin{equation} \label{Eq-cocycle}
\beta(g^{-1}) - \beta(h^{-1}) = -\alpha_{g^{-1}}(\beta(gh^{-1})) = \alpha_{h^{-1}}(\beta(hg^{-1})) \quad \mbox{for} \quad g,h \in \G.
\end{equation}
Let $x = \sum_{g,h} x_{g,h} \otimes e_{g,h}$ be a simple tensor in $L_p(\M \bar\otimes \B(\ell_2(\G)))$. For convenience in the statement and to match the formulation in \cite{JMP}, we assume that $x_{g,h}=0$ whenever $\beta(gh^{-1})=0$, which is equivalent to $\beta(hg^{-1}) = 0$. Then, we deduce the following estimate for $2 \leq p<\infty$ directly from Theorem A with $u_g= \beta(g^{-1})$ together with identities \eqref{Eq-cocycle} and Lemma \ref{colmathom}
\begin{eqnarray*}
\lefteqn{\Big\| \sum_{g,h \in \G} x_{g,h} \otimes e_{g,h} \Big\|_p \approx_{c_p} \Big\| \sum_{g,h \in \G} x_{g,h} \otimes e_{g,h}\otimes \frac{\alpha_{g^{-1}}(\beta(gh^{-1}))}{\|\beta(gh^{-1})\|} \Big\|_{RC_p}} \\ \!\!\!\! & \approx_{c_p}  & \!\!\!\! \Big\| \sum_{g,h \in \G} x_{g,h} \otimes e_{g,h}\otimes \frac{\beta(hg^{-1})}{\|\beta(hg^{-1})\|} \Big\|_{R_p} \! + \Big\| \sum_{g,h \in \G} x_{g,h}\otimes e_{g,h} \otimes \frac{\beta(gh^{-1})}{\|\beta(gh^{-1})\|} \Big\|_{C_p}. 
\end{eqnarray*}                                                               
When $\G$ is amenable and using \cite[Theorem 1.19]{LdlS} and Fourier-Schur transference \ref{Eq-FS}, one recovers that for $x=\sum x_g \otimes \lambda(g) \in L_p(\M \bar\otimes \mathcal{L}(\G))$ with $x_g=0$ when $\beta(g)=0$ and any $2 < p < \infty$
\begin{equation} \label{Eq-JEMS}
\|x\|_p \ \approx_{c_p} \ \Big\| \sum_{g \in \G} x_g \otimes \lambda(g) \otimes \frac{\beta(g)}{\|\beta(g)\|} \Big\|_{C_p} + \Big\| \sum_{g \in \G} x_g \otimes \lambda(g) \otimes \frac{\beta(g^{-1})}{\|\beta(g^{-1})\|}\Big\|_{R_p}.
\end{equation}
This is the main result in \cite{JMP}|the case $1 < p < 2$ follows by duality. This includes nonunimodular groups, not considered in \cite{JMP}. In addition, using local transference from \cite{PRS} we get local forms of \eqref{Eq-JEMS} for nonamenable groups.  

\section{\bf \large Schur multipliers in Schatten $p$-classes}

Given an index set $\Gamma$, the Schur multiplier with symbol $M: \Gamma \times \Gamma \to \C$ is denoted by $S_M$.
The cb-norm on $S_p$ is denoted by $\|S_M\|_{p,\mathrm{cb}}$. Given $1 \le p \le \infty$, we recall the following well known facts
\begin{equation} \label{Schurprop}
\|S_M\|_{p,\mathrm{cb}}=\|S_M\|_{p',\mathrm{cb}}=\|S_{\bar M}\|_{p,\mathrm{cb}}=\|S_{M_{\mathrm{op}}}\|_{p,\mathrm{cb}},
\end{equation}
where $M_{\mathrm{op}}(j,k)=M(k,j)$. Indeed, the equality $\|S_M\|_{p,\mathrm{cb}}=\|S_{M_{\mathrm{op}}}\|_{p',\mathrm{cb}}$ follows by duality. Considering $(\M_{\mathrm{op}},\tau)$ instead of $(\M,\tau)$ gives $\|S_M\|_{p,\mathrm{cb}}=\|S_{M_{\mathrm{op}}}\|_{p,\mathrm{cb}}$ and the last identity comes from a conjugation with the involution. 

\subsection{New criteria for $S_p$-multipliers}

We start with our first criterion.

\demB
Assume first that $u'_k=-u_k$ and $w'_k=-w_k$. We are using Theorem A by fixing some $\H_\R \subset \H$ with associated involution $J$. Let us consider the following maps 
$$\hskip-1.5pt \mathcal{R}_p: S_p(\Gamma) \ni e_{jk} \mapsto e_{jk} \otimes \frac{u_j-u_k}{\|u_j-u_k\|} \in RC_p(\B(\ell_2(\Gamma))),$$
$$\overline{\mathcal{R}}_{p'}: S_{p'}(\Gamma) \ni e_{jk} \mapsto e_{jk} \otimes \frac{J(w_j-w_k)}{\|w_j-w_k\|} \in RC_p(\B(\ell_2(\Gamma))),$$ 
Let $\mathbf{\Lambda}= \mathrm{id} \otimes \Lambda$ be the complete contraction given by homogeneity of $RC_p$. We have for $x\in S_p$ and $y\in S_q$
$${\rm tr} \big( S_M(x)y \big) = \sum_{j,k \in \Gamma} x_{jk} y_{kj} \Big\la \frac{u_j-u_k}{\|u_j-u_k\|},\Lambda \Big(\frac{w_j-w_k}{\|w_j-w_k\|} \Big) \Big\ra = - \big\la \mathbf{\Lambda}^* \mathcal{R}_p(x), \overline{\mathcal{R}}_{p'}(y) \big\ra.$$ Thus $S_M=-\overline{\mathcal{R}}_{p'}^* \mathbf{\Lambda}^* \mathcal{R}_p$ and the estimate follows from Theorem A. The general case follows by a standard $2\times 2$ trick.  Take as index set two disjoints copies $\Gamma_1\cup \Gamma_2$ of $\Gamma$ and set $$\widehat{M}(j,k) = \Big\la \frac{h_j-h_k}{\|h_j-h_k\|}, \Lambda \Big( \frac{g_j-g_k}{\|g_j-g_k\|} \Big) \Big\ra \quad \mbox{with} \quad (h_j,g_j) = \begin{cases} +(u_j,w_j) & \mbox{if} \ j \in \Gamma_1, \\ -(u_j',w'_j) & \mbox{if} \ j \in \Gamma_2. \end{cases}$$ Then, the $\Gamma_1 \times \Gamma_2$ corner of $\widehat{M}$ equals $(M(j,k))_{j,k \in \Gamma}$ and the result follows. \fin

\begin{remark}
\emph{Actually, one can provide a different proof of Theorem B1 i) in the Grothendieck inequality spirit avoiding the Khintchine inequalities but with a worst constant.  We sketch it only for real Hilbert spaces and $\Lambda = \mathrm{id}$.   \cite[Lemma 5.5]{PisSim} gives that $\int_\Omega {\rm sgn}(W(u)){\rm sgn}(W(w))d\mu=\frac 2 \pi \arcsin \la u,w\ra$ for gaussian variables with $\|u\|=\|w\|=1$. Then using Lemma \ref{main0} with a duality argument as above, one gets that the Schur multiplier with symbol $$\arcsin \Big\la \frac{u_j-u_k}{\|u_j-u_k\|},\frac{w_j-w_k}{\|w_j-w_k\|} \Big\ra$$ has norm \hskip-1pt $\lesssim \hskip-1pt \max \{ p, p/p-1 \}^{2}$. We \hskip-1pt can \hskip-1pt remove \hskip-1pt the \hskip-1pt $\arcsin$ \hskip-1pt as \hskip-1pt in \hskip-1pt Krivine's \hskip-1pt proof \hskip-1pt of the Grothendieck inequality \cite[Theorem 5.6]{PisSim} at the price of taking $\sinh$ of the constant.}
\end{remark}

As we have pointed in the Introduction, Grothendieck-Haagerup symbols arise from Theorem B1 i) by fixing $u_k'=0=w_j$. Of course, the constants remain uniformly bounded in that case. In the following result, we show that we can actually improve the constants when only one term vanishes in Theorem B1 i).

\begin{corollary}\label{HGgen2}
The symbols on $\Gamma \times \Gamma$ $$M_1(j,k) = \Big\la \frac{u_j-u_k}{\|u_j-u_k\|},\frac{w_j}{\|w_j\|} \Big\ra \quad \mbox{and} \quad M_2(j,k) = \Big\la \frac{u_j-u_k}{\|u_j-u_k\|},\frac{w_k}{\|w_k\|} \Big\ra$$ satisfy the estimate $\|S_{M_1}\|_{p,\mathrm{cb}} , \|S_{M_2}\|_{p,\mathrm{cb}} \lesssim \max \big\{ p , \frac{p}{p-1} \big\}$ for every $1<p<\infty$. 
\end{corollary}

\dem By \eqref{Schurprop}, it suffices to prove the statement for $M_1$ and $p>2$. We use again the map $T_p$ from the proof of Theorem B1 i). By \eqref{EqRiesz1}, it extends to a cb-map   $S_p(\Gamma) \to C_p(\B(\ell_2(\Gamma)))$ with cb-norm $\lesssim p$. Let $(\delta_j)_{j \in \Gamma}$ denote the canonical basis of $\ell_2(\Gamma)$. It is clear that the map $\rho: S_p(\Gamma) \ni e_{j,k} \mapsto e_{j,k} \otimes \delta_j \in L_p(\B(\ell_2(\Gamma));\ell_2(\Gamma)_c)$ is a complete isometry for all $p \ge 1$. Composing it with $T_p$ and using the associativity of column Hilbert spaces, we get that the map $$V: e_{jk} \mapsto e_{jk} \otimes \Big(\delta_j \otimes \frac{u_j-u_k}{\|u_j-u_k\|} \Big)$$ extends to a cb-map $S_p(\Gamma) \to L_p(\B(\ell_2(\Gamma)); (\ell_2(\Gamma) \otimes_2 \H)_c)$. In addition, the map $U: \ell_2(\Gamma) \otimes_2 \H \to \ell_2(\Gamma)$ given by $U(\delta_j \otimes h) = \la \frac{J(w_j)}{\|w_j\|}, h\ra \delta_j$ is a contraction. By the homogeneity of $C_p$, we get that $(\mathrm{id} \otimes U)V$ is completely bounded, as $J$ is isometric $$(\mathrm{id} \otimes U)V(e_{jk}) = \Big\la \frac{u_j-u_k}{\|u_j-e_k\|},\frac{w_j}{\|w_j\|} \Big\ra e_{jk} \otimes \delta_j = \rho \circ S_{M_1}(e_{jk}).$$ Since $\rho$ is a complete isometry, we conclude $\|S_{M_1}\|_{p,\mathrm{cb}} \le \|(\mathrm{id} \otimes U)V\|_{p,\mathrm{cb}} \lesssim p$. \fin

\begin{remark} \label{RemGralTriang}
\emph{An illustration of Corollary \ref{HGgen2} is given by generalized triangular truncations on $\ell_2(\Z)$.  We consider the family of Schur multipliers whose symbols are given by $\{M_a(j,k) = \mathrm{sgn} (j-a_k): a: \Z\to \Z\}$. This family is $OC_p$-bounded with constant $\max\{p,p/p-1\}$. Indeed, taking $\H_\R=\R$ and $(u_j,u_j',w_j) = (j,-a_j,1)$ the symbols $M_a$ are represented by $$M_a(j,k) = \Big\la \frac {u_j+u_k'}{\|u_j+u_k'\|}, \frac{w_j}{\|w_j\|} \Big\ra = \mathrm{sgn}(j - a_k).$$ Thus, using the $2 \times 2$ trick in the proof of Theorem B1 i) and Corollary \ref{HGgen2}, we get $S_p$-boundedness with constant $C_p=\max\{p,p/p-1\}$. Moreover, arguing as in Corollary \ref{HGgen2} or Remark \ref{examW} below, it follows that this is an $OC_p$-bounded family of operators. Applying this, one gets that the set of Schur multipliers on $\Z$ satisfying $\sup_j \sum_k |M(j,k)-M(j,k+1)| \le 1$ is $OC_p$-bounded by $C_p$. Similarly, multipliers with bounded variation along columns are $OR_p$-bounded with constant $C_p$.}
\end{remark}

\demBB We claim that the symbols 
\begin{equation} \label{Eq-HilbertDividedDifferences}
M(j,k) = \frac{\|\Lambda(w_j) - \Lambda(w_k)\|}{\|w_j - w_k\|}
\end{equation}
define cb-bounded multipliers with $\| S_M \|_{\mathrm{cb}(S_p(\Gamma))} \lesssim \max \big\{p,\frac{p}{p-1} \big\}^{\frac52}$. Indeed, we just apply Theorem B1 i) with $(u_j,u_k',w_j,w_k') = (\Lambda(w_j), -\Lambda(w_k),w_j,-w_k)$. Next we take $\Gamma = \R$ and use \cite[Theorem 1.19]{LdlS} as usual. Consider $\H = L_2(\R)$ and the vectors $w_x = 1_{[0,x]}$ for $x \in \R$. Define $\Lambda: L_2(\R) \to L_2(\R)$ as the pointwise product by $\sqrt{f'}$, which  exists a.e and is bounded by $\sqrt{\|f\|_{\mathrm{lip}}}$. Then $$\frac{\|\Lambda(w_x)-\Lambda(w_y)\|}{\|w_x-w_y\|} = \sqrt{\frac{f(x)-f(y)}{x-y}}.$$ Therefore, the second assertion follows from the first one by approximation. \fin

\begin{remark} \label{Remalphabeta} \emph{In fact, Theorem B1 ii) may be improved as follows. Given any $0 < \beta < 1$ and $f: \R \to \R$ a nondecreasing Lipschitz function we may consider the symbols 
$$M_{f\beta} (x,y) = \Big( \frac{f(x)-f(y)}{x-y} \Big)^\beta.$$ Then, it can be shown for $1 < p < \infty$ 
\begin{equation} \label{Eq-beta1}
\Big\| S_{M_{f\beta}}: S_p(\R) \to S_p(\R) \Big\|_{\mathrm{cb}} \lesssim \max \Big\{ p, \frac{p}{p-1} \Big\}^{2 -\beta} \|f\|_{\mathrm{Lip}}^\beta.
\end{equation}
This follows by the 3-lines lemma with values in $L_p$ and Theorem B2 ii) below. Indeed, let us take $M_{f z}$ for $0 \le \mathrm{Re}(z) \le 1$ and consider the operator-valued analytic function $S_{M_{fz}}(A)$ for some $A \in S_p(\R)$. When $\mathrm{Re}(z) = 0$ it suffices to estimate imaginary powers of $|f(x)-f(y)|$ and $|x-y|$ as Schur multipliers. By standard well-know facts |see for instance \cite[Lemma 1.1]{CGPT2}| the former is dominated by the latter and we get the square of the $S_p$-bound for the Schur multiplier with symbol $|x-y|^{\pm \ii s}$. According to Theorem B2 ii) this gives $(1+|s|^2) \max\{p,p/(p-1)\}^2$ in view of the best constants for this result \cite{CGPT1}. This also follows by Fourier-Schur transference from the classical Mikhlin criterion for completely bounded Fourier $L_p$-multipliers. On the other hand, when $\mathrm{Re}(z) = 1$ we just apply Theorem B2 ii) to obtain the bound $(1+|s|) \max\{p,p/(p-1)\} \|f\|_{\mathrm{Lip}}$, and \eqref{Eq-beta1} follows since the dependance in $s$ is subexponential. Moreover, a similar argument applies for $\alpha$-divided differences of $\alpha$-H\"older functions from \cite[Corollary 3.4]{CGPT1}. Namely, if $0 < \alpha < 1$ and $f: \R \to \R$ is a nondecreasing $\alpha$-H\"older function the symbols $$M_{f\alpha\beta} (x,y) = \frac{|f(x)-f(y)|^\beta}{|x-y|^{\alpha \beta}}$$ define completely $S_p$-bounded multipliers for $|\frac1p - \frac12| < \min \{\alpha, \frac12 \}$.} 
\end{remark}

\begin{remark} \emph{Although \eqref{Eq-beta1} improves Theorem B1 ii), it illustrates how Theorem B1 i) or \eqref{Eq-HilbertDividedDifferences} may be understood as a generalized Hilbert-space-valued divided difference, and other generalizations besides those in Remark \ref{Remalphabeta} are conceivable.}
\end{remark}

\subsection{H\"ormander-Mikhlin-Schur multipliers} \label{S3.2}

We now give a simple proof of the main results in \cite{CLM,CGPT1} from Theorem A. The approach towards HMS multipliers from \cite{CGPT1} is based an a new proof of the H\"ormander-Mikhlin theorem from \cite{JMP,PRS} which relies on Riesz transforms for fractional laplacians. Given $0<\varepsilon<1$, the Hilbert space $W^2_{n,\varepsilon}(\R^n)$ is defined as the completion of $\mathcal{C}_c^\infty(\R^n\setminus\{0\})$ for the norm 
$$\|f\|_{W^2_{n,\varepsilon}} = \big\| \, |\hskip-1pt \cdot \hskip-1pt|^{\frac n2+\varepsilon}\widehat{|\hskip-1pt \cdot \hskip-1pt|^{\varepsilon} f} \, \big\|_{L_2(\R^n)}.$$ 
It consists of bounded continuous  functions and has the property to have a dilation invariant norm. Let us first state Lemma 2.2 from \cite{PRS} in a convenient way for us.

\begin{lemma}\label{convrep}
There exist a complex Hilbert space $\H$ with a strongly continuous unitary action $\alpha: \R^n \curvearrowright \H$ and an $\alpha$-cocycle $\beta: \R^n \to \H$ such that for any $f \in W^2_{n,\varepsilon}(\R^n)$ there exists $u_f \in \H$ such that $$f(x) = \Big\la \frac {\beta(x)}{\|\beta(x)\|}, u_f \Big\ra \quad \mbox{and} \quad \|u_f\|_\H \approx_{n,\varepsilon} \|f\|_{W^2_{n,\varepsilon}}.$$ In particular, we find $f(x-y) = \big\la \frac {\beta(x) - \beta(y)}{\|\beta(x)-\beta(y)\|}, \alpha_{-y} (u_f) \big\ra = \big\la \frac {\beta(x) - \beta(y)}{\|\beta(x) - \beta(y)\|}, -\alpha_{-x} (u_f) \big\ra$.
\end{lemma}

\begin{remark} \label{examW}
\emph{Therefore, elements in $W^2_{n,\varepsilon}(\R^n)$ define Toeplitz Schur multipliers thanks to Corollary \ref{HGgen2}.  Moreover, let $(f_x)_{x\in \R^n}$ be a family of functions in the unit ball of $W^2_{n,\varepsilon}(\R^n)$ for some $0<\varepsilon<1$. Then, the Schur multiplier on $\R^n$ with symbol $$M(x,y) = f_x(x-y)$$ for $x,y \in \R^n$ satisfies $\|S_M\|_{p,\mathrm{cb}} \lesssim_{n,\varepsilon} C_p$ for $1<p<\infty$. In fact, the family of such Schur multipliers is $OC_p$-bounded with the same constant since a map satisfying $T(x\otimes \delta_\ell) = S_{M_\ell}(x)\otimes \delta_n$ corresponds to the restriction of a Schur multiplier $S_M$ on $\B(\ell_2(\Gamma) \otimes_2 \ell_2(\Z))$. Indeed, assume that $M_\ell$ is given by a sequence $(w^\ell_j)_{j \in \Gamma}$, then $M$ is given by \vskip-15pt $$M((j,\ell)(k,m)) = \Big\la \frac{u_j-u_k}{\|u_i-u_k\|}, \frac{w_j^\ell}{\|w_j^\ell\|} \Big\ra.$$}
\end{remark}

\begin{remark} \label{Rem-Sobolev}
\emph{According to  \cite[Lemma 2.3]{PRS} we know that, for functions with support in a compact set not containing $0$, the $W^2_{n,\varepsilon}$-norm is smaller than the Sobolev norm $$\|f\|_{H_{\frac n 2 +\varepsilon}} = \big\| \, (1+|\hskip-1pt \cdot \hskip-1pt |^2)^{\frac n4 + \frac \varepsilon 2} \widehat{f} \ \big\|_{L_2(\R^n)}$$ and compactly supported functions in $H_{\frac d 2 +\varepsilon}$ define Toeplitz Schur multipliers. This is obvious though, since the Fourier transform of elements in this Sobolev space belong to $L_1(\R^n)$. In fact, if $(f_x)$ are functions in the unit ball of $H_{\frac d 2+\varepsilon}$ we get a $OC_\infty$-bounded family without any assumption on the support because $$f_x(x-y)= \int_{\R^n}   (1+|s|^2)^{\frac d 4+ \frac \varepsilon 2}\widehat{f_x}(s)e^{2\ii \pi \la s,x\ra} \frac 1{(1+|s|^2)^{\frac d 4+ \frac \varepsilon 2}} e^{-2\ii \pi \la s,y\ra} ds = \big\la \xi_x,\eta_y \big\ra_{L_2},$$ 
where $\sup_{x}\|\xi_x\|_2\lesssim_{d,\varepsilon} 1$ and $\sup_{y}\|\eta_y\|_2\lesssim_{d,\varepsilon} 1$. Thus the $OC_\infty$-boundedness comes from the Grothendieck-Haagerup criterion (as $\eta_y$ does not depend of $f_x$) and we also get $OC_p$-boundedness with a constant independent of $p$ for every $p \ge 1$. This illustrates that the spaces $W_{n, \varepsilon}^2(\R^n)$ are definitely larger, allowing more singularities on the symbols considered. Therefore a H\"ormander-Mikhlin theorem in terms of these spaces |like the one given below| is formally stronger.}
\end{remark}

In what follows and again by approximation, we shall work with $\Gamma = \R^n$. Take $\Delta_k = \{(x,y) \in \R^n \times \R^n: 2^k \le |x-y|_\infty < 2^{k+1}\}$. This is a Littlewood-Paley decomposition which arises as a combination of Hilbert transforms (one-dimensional Riesz transforms) that satisfies \eqref{LP} with $D_p \lesssim \max \{p,p/p-1\}^{2n}$. Other smoother choices of Littlewood-Paley decomposition yield better constants, but we shall not consider them for the sake of simplicity on our presentation. In the following result we fix $0<\varepsilon<1$ and $\phi \in \mathcal{C}_c^\infty(\R^n \setminus \{0\})$ even and identically 1 on $\Delta_0$.

\begin{theorem} \label{Thm-StrongHMS}
Let $M:\R^n \times \R^n \to \C$ be a continuous function such that 
$$\sup_{x \in \R^n} \sup_{k \in \Z} \big\| \phi(2^k \cdot) M(x, x-\cdot) \big\|_{W^2_{n,\varepsilon}} + \sup_{y\in \R^n} \sup_{k \in \Z} \big\| \phi(2^k \cdot) M(y - \cdot,y) \big\|_{W^2_{n,\varepsilon}} \le 1.$$
Then $S_M$ defines a completely bounded Schur multiplier on $S_p$ for every $1<p<\infty$.
\end{theorem}

\dem Let $M_k(x,y)= \phi (2^{-k}(x-y)) M(x,y)$ for $x,y \in \R^n$ and $k \in \Z$. It is clear that $M(x,y) = \sum_{k \in \Z} 1_{(x,y) \in \Delta_k} M_k(x,y)$ for $x \neq y$ because $\varphi(2^{-k}(x-y)) = 1$ if $(x,y) \in \Delta_k$. Consider the functions $$f_{x,k}(y)= \phi(2^{-k}y)M(x,x-y) \quad \mbox{for} \quad (x,k) \in \R^n \times \Z.$$ We have $M_k(x,y) = f_{x,k}(x-y)$ and our assumption means that $\|f_{x,k}\|_{W^2_{n,\varepsilon}} \le 1$. In particular, the family $\{S_{M_k}: k \in \Z\}$ is $OC_p$-bounded by Remark \ref{examW}. Similarly let $g_{k,y}(x) = \phi(2^{-k}x) M(y-x,y)$. Then $M_k(x,y) = g_{k,y}(y-x)$ and the same argument shows that indeed the family $\{S_{M_k}: k \in \Z\}$ is $OR_p$-bounded, and thus $ORC_p$-bounded. The assertion then follows from  Lemma \ref{RCbound}. \fin

\demBBBB The same argument with the generalized triangular projections from Remark \ref{RemGralTriang} gives the Marcinkiewicz type criterion of \cite{CLM}. Indeed, any multiplier of strong 1-variation along rows sits in the closure of the absolute convex hull of generalized triangular projections and the identity. Thus one can use Lemma \ref{absconv} for the $OC_p$-boundedness and similarly for $OR_p$-boundedness. This argument is very close to the original one given in \cite{CLM}. \fin

\demBBB By \cite[Lemma 2.2]{PRS}, the criterion in Theorem \ref{Thm-StrongHMS} is stronger than the main result in \cite{CGPT1} in both the H\"ormander-Mikhlin and Sobolev formulations. This already justifies Theorem B2 ii). Our argument yields worse constants but,  taking a smooth Littlewood-Paley partition and Remark \ref{Rem-Sobolev} into account, we recover the optimal behavior of constants $C_p = \max \{p, p/p-1\}$. \fin

\begin{remark}
\emph{When $\Gamma \neq \R^n$ is another locally compact group, the results around Lemma \ref{convrep} remain true in some sense \cite{JMP}. Thus, the technique in this paper is also valid to recover \cite[Corollary 4.5]{CGPT2} for nonorthogonal actions, which in turn generalizes the results from \cite{JMP0,JMP}.}
\end{remark}

\noindent \textbf{Acknowledgement.} We would like to thank Mikael de la Salle, Fedor Sukochev and Dimitry Zanin for helpful comments which have improved the presentation.

\bibliographystyle{amsplain}

\begin{thebibliography}{99}

\bibitem {A} J. Arazy, Schur-Hadamard multipliers in the spaces $C_p$. Proc. Amer. Math. Soc. \textbf{86} (1982), 59-64.

\bibitem {AK} C. Arhancet and C. Kriegler, Riesz transforms, Hodge-Dirac operators and functional calculus for
multipliers. Lecture Notes in Mathematics \textbf{2304}. Springer, 2022.

\bibitem {B1} D. Bakry, Transformation de Riesz pour les semi-groupes sym\'etriques. S\'eminaire de Probabilit\'es XIX. Lecture Notes in Math. \textbf{1123} (1985), 130-175.

\bibitem {B2} D. Bakry, Etude des transformations de Riesz dans les vari\'et\'es Riemanniennes \`a courbure de
Ricci minor\'ees. S\'eminaire de Probabilit\'es XXI. Lecture Notes in Math. \textbf{1247} (1987), 137-172.

\bibitem {DCH} J. de Canni\`ere and U. Haagerup, Multipliers of the Fourier algebras of some simple Lie groups and their discrete subgroups. Amer. J. Math. \textbf{107} (1985), 455-500.

\bibitem{CdlS} M. Caspers and M. de la Salle, Schur and Fourier multipliers of an amenable group acting on non-commutative $L_p$-spaces. Trans. Amer. Math. Soc. \textbf{367} (2015), 6997-7013.

\bibitem {CLM} C.Y. Chuah, Z.C. Liu and T. Mei, A Marcinkiewicz multiplier theory for Schur multipliers. Anal. \& PDE, to appear.

\bibitem {CGPT1} J.M. Conde-Alonso, A.M. Gonz\'alez-P\'erez, J. Parcet and E. Tablate, Schur multipliers in
Schatten–von Neumann classes. Ann. of Math. \textbf{198} (2023), 1229-1260.

\bibitem {CGPT2} J.M. Conde-Alonso, A.M. Gonz\'alez-P\'erez, J. Parcet and E. Tablate, A Hörmander-Mikhlin theorem for high rank simple Lie groups. Preprint 2022, arXiv:2201.08740.

\bibitem {CH} M. Cowling and U. Haagerup, Completely bounded multipliers of the Fourier algebra of a simple Lie group of real rank one. Invent. Math. \textbf{96} (1989), 507-549.


\bibitem {G} A. Grothendieck, R\'esum\'e de la th\'eorie m\'etrique des produits tensoriels topologiques. Boll. Soc. Mat. Sao-Paulo \textbf{8} (1956), 1-79.

\bibitem {Gu} R. Gundy, Sur les transformations de Riesz pour le semi-groupe d’Ornstein-Uhlenbeck. C.R. Acad. Sci. Paris \textbf{303} (1986), 967-970.

\bibitem {GV} R. Gundy and N. Varopoulos, Les transformations de Riesz et les int\'egrales stochastiques. C.R. Acad. Sci. Paris \textbf{289} (1979), 13-16.

\bibitem {H} U. Haagerup, An example of a non nuclear $C^*$-algebra, which has the metric approximation property. Invent. Math. \textbf{50} (1979), 279-293.


\bibitem {HvNVW} T. Hyt\"onen, J. van Neerven, M. Veraar and L. Weis, Analysis in Banach Spaces | Volume II: Probabilistic Methods and Operator Theory. Springer, 2017. 

\bibitem{JLMX} M. Junge, C. Le Merdy and Q. Xu, $H^\infty$ functional
  calculus and square functions on noncommutative $L^p$-spaces. {Ast\'erisque}. \textbf{305} (2006).
  
\bibitem {JM} M. Junge and T. Mei, Noncommutative Riesz transforms—A probabilistic approach. Amer. J. Math. \textbf{132} (2010), 611-681.

\bibitem {JMP0} M. Junge, T. Mei and J. Parcet, Smooth Fourier multipliers in group von Neumann algebras. Geom. Funct. Anal. \textbf{24} (2014), 1913-1980..

\bibitem {JMP} M. Junge, T. Mei and J. Parcet, Noncommutative Riesz transforms—dimension free bounds and
Fourier multipliers. J. Eur. Math. Soc. \textbf{20} (2018), 529-595.

\bibitem {dLdlS} T. de Laat and M. de la Salle, Approximation properties for noncommutative $L_p$ of high rank lattices and nonembeddability of expanders. J. Reine Angew. Math. \textbf{737} (2018), 46-69.

\bibitem {LdlS} V. Lafforgue and M. de la Salle, Noncommutative $L_p$-spaces without the completely bounded approximation property. Duke. Math. J. \textbf{160} (2011), 71-116.

\bibitem {LP0} F. Lust-Piquard, In\'egalit\'es de Khintchine dans $C_p$ ($1 < p < \infty$). C.R. Acad. Sci. Paris \textbf{303} (1986), 289-292.

\bibitem {LP1} F. Lust-Piquard, Riesz transforms associated with the number operator on the Walsh system
and the fermions. J. Funct. Anal. \textbf{155} (1998), 263-285.

\bibitem {LP2} F. Lust-Piquard, Riesz transforms on deformed Fock spaces. Comm. Math. Phys. \textbf{205} (1999),
519-549.

\bibitem {LP3} F. Lust-Piquard, Dimension free estimates for discrete Riesz transforms on products of abelian
groups. Adv. Math. \textbf{185} (2004), 289-327.

\bibitem {LPP} F. Lust-Piquard and G. Pisier, Non-commutative Khintchine and Paley inequalities. Ark. Mat. \textbf{29} (1991), 241-260.

\bibitem {M} P.A. Meyer, Transformations de Riesz pour les lois gaussiennes. S\'eminaire de Probabilit\'es XVIII. Lecture Notes in Math. \textbf{1059} (1984), 179-193.

\bibitem {N} A. Naor, Discrete Riesz transforms and sharp metric $X_p$ inequalities. Ann. of Math. \textbf{184} (2016), 991-1016, 2016.

\bibitem {NR} S. Neuwirth and \'E. Ricard, Transfer of Fourier multipliers into Schur multipliers and sumsets in a discrete
group. Canad. J. Math. \textbf{63} (2011), 1161-1187.

\bibitem {PRS} J. Parcet, \'E. Ricard and M. de la Salle, Fourier multipliers in $SL_n(\mathbf{R})$. Duke Math. J. \textbf{171} (2022), 1235-1297.

\bibitem {PST} J. Parcet, M. de la Salle and E. Tablate, The local geometry of idempotent Schur multipliers. Preprint 2024, arXiv 2312.02895.

\bibitem {P0} G. Pisier, Riesz transforms: a simpler analytic proof of P.A. Meyers inequality. In S\'eminaire de Probabilit\'es XXII. Lecture Notes in Math. \textbf{1321} (1988), 485501.

\bibitem {PisSim} G. Pisier, Similarity Problems and Completely Bounded Maps. Lecture Notes in Mathematics \textbf{1618}. Springer-Verlag, 2001.

\bibitem {P} G. Pisier, Introduction to Operator Space Theory. London Mathematical Society Lecture Note Series \textbf{294}. Cambridge University Press, 2003.

\bibitem {PisBAMS} G. Pisier, Grothendieck’s Theorem, past and present. Bull. Amer. Math. Soc. \textbf{49} (2012), 237-323. 

\bibitem {PS} D. Potapov and F. Sukochev, Operator-Lipschitz functions in Schatten-von Neumann classes. Acta Math. \textbf{207} (2011), 375-389.

\bibitem {S} E.M. Stein, Some results in harmonic analysis in $\R^n$, for $n \to \infty$. Bull. Amer. Math. Soc. \textbf{9} (1983), 71-73.
\end{thebibliography}

\enlargethispage{2,5cm}

\vskip3pt

\small

\noindent \textbf{Adri\'an M. Gonz\'alez-P\'erez} \hfill \textbf{Javier Parcet} \\ ICMAT \hfill ICMAT \\ Universidad Aut\'onoma de Madrid \hfill Consejo Superior de Investigaciones Cient{\'\i}ficas 
\\ \texttt{adrian.gonzalez@uam.es} \hfill \texttt{parcet@icmat.es}

\noindent \textbf{Jorge P\'erez-Garc\'ia} \hfill \noindent \textbf{\'Eric Ricard} \\
ICMAT \hfill LMNO \\ Consejo Superior de Investigaciones Cient{\'\i}ficas \hfill Universit\'e de Caen Normandie  
\\ \texttt{jorge.perez@icmat.es} \hfill \texttt{eric.ricard@unicaen.fr}




\end{document}